\documentclass{elsarticle}
\usepackage{amsmath,amssymb,calc,enumerate,booktabs}
\usepackage{amsthm}
\usepackage{algorithmic}
\usepackage{algorithm}
\usepackage{units}
\usepackage{tikz}
\usepackage{pgfplots}
\usepackage{caption}
\usepackage{subcaption}
\usepackage{rotating}

\pgfplotsset{compat=newest}

\usetikzlibrary{positioning,shapes,shadows,arrows}

\title{Satellite downlink scheduling problem: A case study\tnoteref{t1}}

\tnotetext[t1]{NOTICE: this is the author's version of a work that was accepted for publication in Omega.  Changes resulting from the publishing process, such as peer review, editing, corrections, structural formatting, and other quality control mechanisms may not be reflected in this document.  Changes may have been made to this work since it was submitted for publication.  A definitive version was subsequently published in Omega 53, 2015, DOI 10.1016/j.omega.2015.01.001.}

\author{Daniel~Karapetyan\fnref{not,sfu}\corref{cor1}}
\ead{daniel.karapetyan@gmail.com}

\author{Snezana~Mitrovic~Minic\fnref{sfu,mda}}

\author{Krishna~T.~Malladi\fnref{sfu}}

\author{Abraham~P.~Punnen\fnref{sfu}}

\fntext[not]{ASAP Research Group, School of Computer Science, University of Nottingham, Nottingham NG8\,1BB, UK}
\fntext[sfu]{Department of Mathematics, Simon Fraser University, Surrey, BC V3T\,0A3, Canada}
\fntext[mda]{MDA Systems Ltd., Richmond, BC V6V\,2J3, Canada}

\cortext[cor1]{Corresponding author}

\date{}

\renewcommand{\comma}{,\allowbreak\ }

\begin{document}

\begin{abstract}
The synthetic aperture radar (SAR) technology enables satellites to efficiently acquire high quality images of the Earth surface.  This generates significant communication traffic from the satellite to the ground stations, and, thus, image downlinking often becomes the bottleneck in the efficiency of the whole system.  In this paper we address the downlink scheduling problem for  Canada's Earth observing SAR satellite, RADARSAT-2.  Being an applied problem, downlink scheduling is characterised with a number of constraints that make it difficult not only to optimise the schedule but even to produce a feasible solution.  We propose a fast schedule generation procedure that abstracts the problem specific constraints and provides a simple interface to optimisation algorithms.  By comparing empirically several standard meta-heuristics applied to the problem, we select the most suitable one and show that it is clearly superior to the approach currently in use.
\end{abstract}


\maketitle

\section{Introduction}

Efficient scheduling of image acquisition and image downlinking plays a vital role in satellite mission planning.  These operations are often interlinked and solved using scheduling heuristics. Most of the literature on satellite mission planning (image acquisition and downlinking) is divided into two categories: optical satellites~\cite{a10,z5} and Synthetic Aperture Radar (SAR) satellites~\cite{a1,a12,z2,z4}.

This paper deals with the downlink scheduling portion of the mission planning operations of Canada's Earth observing SAR satellite, RADARSAT-2.
We assume that image acquisition schedule is given and not to be changed. This is consistent with existing policies and practices.

Recently, there is a steady increase in demand for RADARSAT-2 imagery.  Thus any improvements in image downlink operations would improve the efficiency of the RADARSAT-2 mission and this is the primary motivation behind this study.  We have addressed only downlink scheduling problem because in this satellite mission the image acquisition scheduling is performed directly by the customers and it is ruled by the customer priorities.  The customers have direct access to an image scene ordering tools, and they select and order particular scenes (at particular time) based on their area of interest.  An order from a high-priority customer is always given priority.  If such an order generates cancellations of previously submitted orders, the affected low-priority customers are alerted to repeat their ordering procedure.  

It is not always possible to schedule downlinks for all the images.  If some image could not be downlinked by the deadline assigned of the order, the satellite operations centre revises the image acquisition plan by removing that image.  Customers affected by such cancellations are alerted to place new orders.

Currently used downlink scheduling process exploits a greedy-like algorithm \cite{greedy-like} followed by human intervention whenever necessary.  We explored local search heuristics and metaheuristics to improve the efficiency of the downlink scheduling.  Our experimental study on real-world problem instances has shown that the proposed techniques significantly improved downlink throughput and schedule quality.

The satellite image downlink scheduling problem and its variations have been studied by many authors.  Some of these works were focused on case studies for specific satellites or space missions~\cite{z4,donati2011,our-chapter} whereas others are more general purpose in nature~\cite{a4,a2,a5,a6,a3,a9,z3,z5,a11}. Literature from machine scheduling~\cite{b1,a13,Sterna2011,Sawik2010} and resource-constrained project scheduling~\cite{rcpsp,rcpsp-application} are also relevant in solving the satellite image downlink scheduling problem (SIDSP)\@.  However, each mission planning problem has its own restrictions and properties that can be exploited.  For example, some studies consider single satellite problems while others deal with satellite constellations.  Special care is needed to make sure that the restrictions are handled adequately which sometimes changes the inherent combinatorial structure of the problem significantly.  Thus, investigating case studies of special SIDSPs are interesting and relevant as established in this study, although existing literature on the SIDSP considerably influenced our work.

The paper is organised as follows. Section~\ref{sec:practicalProblem} describes the real-world problem, and Section~\ref{sec:problem} introduces its mathematical model as well as various notations and definitions.  Sections~\ref{sec:dev_algorithm} and~\ref{sec:ScheduleGenerator} deal with our heuristic algorithms.  Data analysis and test instances are reported in Section~\ref{sec:testbed} followed by computational results in Section~\ref{sec:results} and concluding remarks in Section~\ref{sec:conclusion}.

\section{The RADARSAT-2 SIDSP}
\label{sec:practicalProblem}

We consider the problem arising in satellite industry that deals with scheduling downlinks of images assuming that the schedule for image acquisition is already generated. We study RADARSAT-2 satellite that orbits the Earth to acquire images that can be downlinked to a set $G$ of stationary ground stations for further processing.  Since the problem has a lot in common with scheduling problems and, in particular, with the Resource Constrained Project Scheduling Problem, our notations are close to the ones used in the scheduling literature.

Let $V$ be a set of $n$ downlink requests to be scheduled within the planning horizon of 24 hours.  For each request $j \in V$, release time $r_j$, deadline $d_j$, downlink duration $p_j$, priority $w_j$ and ground station $g_j \in G$ are prescribed.  The interval $[r_j, d_j]$ is called the \emph{time window} of request $j$.  Downlink requests are classified as \emph{regular} and \emph{urgent}.  The downlink of urgent requests has to start as close as possible to their release times.  If several urgent requests are competing, then their priorities have to be taken into account.  However, downlinks of regular requests are more flexible and are primarily governed by their priorities. Urgent requests have absolute priority dominance over the regular (non-urgent) requests, i.e.\ any (small) improvement in  downlinking of urgent requests is preferred over large improvements in  downlinking of regular requests.  Finally, some requests have to be downlinked to two ground stations and we call them \emph{dual requests}.  Each dual request is represented by two requests $i, j \in V$, and we are given a set $D$ of pairs $(i, j)$ of dual requests.

A downlink activity can be carried out only when the satellite is passing over a ground station.  This time interval is called \emph{visibility mask} of the station.  The RADARSAT-2 has two antennas for downlinking, and it can work in half-power and full-power setting.  When working in half-power setting, the two antennas can work separately and they can independently downlink two different images to one or two different ground stations simultaneously.  In full-power setting, the satellite can process only one downlink at a time.

Each ground station $g \in G$ has one or two channels for receiving the downlinked images.  Ground stations can also be classified based on their transmission power.  Let $G_1$ be the set of ground stations that are in \emph{half-power setting} and $G_2$ be the set of ground stations that are in \emph{full-power setting}.  Then $G = G_1 \cup G_2$ and $G_1 \cap G_2 = \emptyset$.  When two images are downlinked one after another to stations with the same power setting, there must be a gap (set up time) $\delta$ time units between the two downlinks.  When two images are downlinked consecutively to two stations with different power settings, the required gap between the downlinks is $\Delta > \delta$.

The visibility masks of a ground station $g \in G$ can be represented by a collection $M_g$ of non-overlapping time intervals called the \emph{normal visibility masks}.  Certain downlink requests require better reliability and they have to be downlinked within \emph{high reliability visibility mask} $M^1_g$.  Each high reliability visibility mask $m^1 \in M^1_g$ is a sub-interval of some normal visibility mask $m \in M_g$.

\section{Mathematical Model}
\label{sec:problem}

Note that SIDSP deals with the problem of finding an image downlink schedule so that a utility function is maximized.  The utility function considers the number of downlinks scheduled, their priority values, and the difference between the downlink start time and the start time of the request time window (tardiness).  Rejection of requests is allowed, and a rejected request is referred to as \emph{unscheduled}.  Note that the downlink scheduling problem is usually oversubscribed due to large customer demand for satellite imagery.  Thus, if request rejection is not allowed, the downlink scheduling problem would often be infeasible.  Furthermore, the request rejection assumption allows the scheduling algorithm to choose the requests that maximise the resource utilisation.  

A solution $S$ to SIDSP --- a \emph{schedule} --- is a set $S \subseteq V$ of scheduled requests and associated downlink start times $S_j$ for each $j \in S$.  Our model for the SIDSP is to:
\begin{align*}
\text{Maximize } & f(S) = \sum_{j \in S} w_j \left( 1 - \alpha \cdot \frac{S_j - r_j}{d_j - p_j - r_j} \right) \\
\text{subject to } & S \in \mathbb{F} \,, \\
\end{align*}
where $0 \le \alpha \le 1$ is a parameter reflecting the importance of tardiness minimisation ($\alpha = 0$ disables tardiness minimisation while $\alpha = 1$ means that scheduling a request to the end of its time window is as bad as rejecting it), and the collection $\mathbb{F}$ contains all schedules $S$ satisfying the following constraints:
\begin{enumerate}[(1)]
\item No downlink activity happens outside the planning horizon.
\item A downlink cannot start earlier than the release time of the corresponding  request and must be finished by its deadline: $r_j \le S_j \le d_j - p_j$ for each $j \in S$.
\item Once a downlink starts, it cannot be preempted, i.e.\ an image cannot be split into several fragments.
\item Each urgent request has to be scheduled at the earliest possible time in its time window even at the cost of delaying or cancelling some regular downlinks.  In several urgent requests are competing, their priorities have to be taken into account.
\item For each dual request $(i, j) \in D$, $i \in S$ if and only if $j \in S$.
\item \label{constraint:delta} There must be a gap of at least $\delta$ units between two consecutive downlinks under the same power setting.
\item There must be a gap of at least $\Delta \ge \delta$ time units between two consecutive downlinks under different power setting.
\item If the satellite is in the full-power setting, then only one antenna can work and the other has to be idle.  For both of the satellite antennas to be working independently, the satellite has to be in the half-power setting.  The full-power setting is used if and only if the satellite transmits data to a full-power ground station.
\end{enumerate}

\noindent We assume the following:
\begin{enumerate}[(1)]
\item The satellite is in the half-power setting at the beginning of the planning horizon and that it has to be in the half-power setting at the end as well.
\item Downlinking (in half-power setting) can  start right from the beginning of the planning horizon, i.e.\ it is guaranteed that no half-power downlink activity happened within $\delta$ time units before the beginning of the planning horizon.
\item We tackle constraint (\ref{constraint:delta}) by the following pre-processing procedure:
Add $\delta$ to $p_j$ and $d_j$ for all requests $j \in V$.  Also, add $\delta$ to the upper bounds of each of the time intervals in normal and high reliability visibility masks for all $g \in G$.  Extend the planning horizon by $\delta$. Further, replace $\Delta$ with $\Delta - \delta$.
\item The solution mechanism ensures that the urgent requests get absolute priority over regular requests.  We discuss our approach to this issue in Section~\ref{sec:dev_algorithm}.
\end{enumerate}


The SIDSP is NP-hard since several NP-hard machine scheduling problems are special cases of it.  Consider, for example, the parallel machine scheduling problem with two machines representing satellite antennas.

Heuristic algorithms for the scheduling problems of this class usually exploit the so-called \emph{serial scheduling scheme}, in which the solver searches in the space of job sequences while a polynomial time schedule generator converts the job sequences into schedules.  The schedule generator is a greedy algorithm scheduling the jobs (requests in our work) in the given order, choosing the earliest available position for each of them.  An important property of the schedule generator is that it always generates active schedules, i.e.\ schedules such that none of the unscheduled jobs can be added to it and no scheduled job can be advanced without delaying some other job~\cite{kolisch}.  Given that our objective function is \emph{regular} (i.e.\ delaying or deleting a downlink cannot improve the solution if no other changes are introduced), there exists a sequence of requests generating an optimal schedule~\cite{project-scheduling}.

Moreover, for any active schedule, there exists a sequence generating that schedule (indeed, it is enough to sort the jobs in ascending order by their start times).  This implies that the schedule generator can produce a worst possible active schedule.  For example, it can produce a schedule $S$ of objective $f(S) = 0$ by scheduling requests of zero priority and leaving the requests with non-zero priorities unscheduled.



With the serial scheduling scheme in mind, the SIDSP can be viewed as a problem of optimising a permutation of the downlink requests for which a schedule will be generated in the following manner: From the order defined by the permutation, schedule each request to the earliest available time considering the constraints related to the satellite antennas, ground station channels, visibility masks and dual requests.  
Let $\wp$ be the set of all permutations of elements of $V$. For any $\pi \in \wp$, let $S_{\pi}$ be the corresponding schedule generated by the schedule generator algorithm (see Section \ref{sec:ScheduleGenerator}).
Thus, to solve SIDSP, we solve the following \emph {Downlink Request Permutation Problem} (DRPP):
\begin{align*}
\text{Maximize } & \phi(\pi) =  f(S_{\pi}) \\
\text{subject to } & \pi \in \wp.
\end{align*}
Note that to evaluate the solution quality of a permutation we need to generate the schedule from that permutation, i.e.\ to apply the schedule generator algorithm described in Section~\ref{sec:ScheduleGenerator}.

\section{Solution Approach}
\label{sec:dev_algorithm}

Instead of dealing with urgent requests separately in the objective function or in some other way, we used a two phase solution approach.
The first phase schedules all urgent requests, and the second phase schedules all regular requests using the remaining resources.  Hence, our approach respects the requirement to give the urgent requests ultimate priority over the regular requests.  A high-level description of the algorithm is given as follows:
\begin{description}
\item[Phase 1 --- Urgent requests:] Schedule all urgent requests using heuristic $\mathcal{H}$.
\item[Phase 2 --- Regular requests:] Fix the urgent requests scheduled in Phase~1, update resource availability accordingly and schedule regular requests using heuristic $\mathcal{H}$.
\end{description}

In the following sections we describe several heuristic algorithms for the DRPP that can be used as  $\mathcal{H}$\@.  We primarily focused on standard algorithmic paradigms since one of the objectives was to propose  algorithms that are easy to understand and implement.  In particular, we limited our experiments to the Greedy Randomised Adaptive Search Procedure (GRASP), Ejection Chain, Simulated Annealing and Tabu Search algorithms.

\subsection{Construction Heuristic}
\label{sec:InitialSoln}

One of the components of our heuristics is a construction algorithm, that generates an initial schedule from the list of sorted requests in $V$ (the so-called \emph{priority-rule based scheduling method} known to be efficient for quick generation of reasonably good solutions~\cite{project-scheduling}).  Despite significant achievements in designing complicated dispatching rules for standard scheduling problems (see, e.g., \cite{dispatching-rules}), we preferred a simple dispatching scheme as the quality of the initial solutions is not crucial to us.  We considered the following sorting criteria and tie breakers:
\begin{enumerate}
	\item Priority $w_j$;
	\item Time window duration $d_j - r_j - p_j$\;
	\item Downlink duration $p_j$.
\end{enumerate}
After extensive computational experiments using various combinations of the listed criteria, we selected sorting the elements of $V$ by $w_j$ in the descending order, ties broken by $\lfloor d_j - r_j - p_j \rfloor$ in the ascending order.




\subsection{GRASP}

The Greedy Randomised Adaptive Search Procedure is a simple meta-heuris\-tic often applied to scheduling problems~\cite{marti}.  GRASP repeatedly generates solutions with a randomised greedy constructor followed by a local search phase.  Due to the randomness of the greedy procedure, GRASP is likely to produce new solutions on every iteration.  The best out of all the produced solutions is selected in the end.  For details, see Algorithm~\ref{alg:grasp}.

\begin{algorithm}[htb]
\caption{GRASP Algorithm}
\label{alg:grasp}

\begin{algorithmic}
\REQUIRE Given time $\mathcal{T}$
\ENSURE A permutation optimised with respect to $\phi()$
\WHILE {$\text{``elapsed time''} < \mathcal{T}$}
	\STATE $\pi' \gets \mathit{GreedyRandomisedConstructor()}$
	\STATE $\pi' \gets \mathit{LocalSearch}(\pi')$
	\IF {$\phi(\pi') > \phi(\pi)$}
		\STATE $\pi \gets \pi'$
	\ENDIF
\ENDWHILE
\RETURN $\pi$
\end{algorithmic}
\end{algorithm}

As a randomised greedy constructor $\mathit{GreedyRandomisedConstructor()}$, we use a modification of the construction procedure described in Section~\ref{sec:InitialSoln}.  In particular, on every iteration of building $\pi$, $\mathit{GreedyRandomisedConstructor()}$ orders all the remaining downlink requests as described in Section~\ref{sec:InitialSoln} and selects one of the first ten candidates randomly with uniform probability distribution.

Our local search $\mathit{LocalSearch}(\pi)$ explores a \emph{swap neighbourhood}.  The swap neighbourhood $N_{swap}(\pi)$ consists of all the solutions that can be obtained from the permutation $\pi$ by swapping two of its elements.  The size of the neighbourhood is $|N_{swap}(\pi)| = \frac{n (n - 1)}{2}$, and, hence, it would take $O(n^5)$ time to explore it (we will show in Section~\ref{sec:ScheduleGenerator} that the complexity of the schedule generator algorithm is $O(n^3)$).  With such a high complexity of the local search, GRASP is likely to perform only few iterations, which is not enough to exploit the strength of the meta-heuristic.  To speed up the local search phase, we decided to explore the neighbourhood in a random order and terminate the search when a prescribed time limit is reached.  For details, see Algorithm~\ref{alg:random-ls}.

\begin{algorithm}[htb]
\caption{GRASP local search $\mathit{LocalSearch}(\pi)$}
\label{alg:random-ls}

\begin{algorithmic}
\REQUIRE initial permutation $\pi^0$; the time given for one run of the local search
\ENSURE Improved permutation with respect to $\phi()$

\STATE $\pi \gets \pi^0$
\WHILE {given time did not elapse}
	\STATE Select randomly $i \neq j \in \pi$ with uniform probability distribution
	\STATE $\pi' \gets \mathit{Swap}(\pi, i, j)$
	\IF {$\phi(\pi') > \phi(\pi)$}
		\STATE $\pi \gets \pi'$
	\ENDIF
\ENDWHILE
\RETURN $\pi$
\end{algorithmic}
\end{algorithm}

Here and in the rest of the paper function $\mathit{Swap}(\pi, i, j)$ swaps elements in positions $i$ and $j$ of the permutation $\pi$; the original permutation $\pi$ remains unchanged.

\subsection{Ejection Chain Algorithm}
\label{sec:ejection_chain}

Ejection chain methods~\cite{a16} have commonly been used in developing Very Large Scale Neighbourhood (VLSN) search algorithms~\cite{a15,a14} for solving complex combinatorial optimisation problems. For example, the well known Lin-Kernighan heuristic --- an efficient heuristic for solving the travelling salesman problem~\cite{Gamboa2006} --- is an ejection chain algorithm.  We use the idea of the ejection chains to develop a simple and effective heuristic to solve the DRPP.

The data structure used in our ejection chain algorithm is a pair $(\pi, h)$, where $\pi$ is a permutation of all the requests in $V$ and $h \in \{ 1, 2, \ldots, n \}$ is a position in this permutation called a \emph{hole}.  In order to calculate the objective value $\xi(\pi, h)$, copy all the permutation $\pi$ skipping the element in position $h$ and feed this copy to the schedule generator algorithm (Section~\ref{sec:ScheduleGenerator}) to obtain the schedule and calculate its objective.



The basic move in our ejection chain algorithm is swapping the element of $\pi$ in position $i$ with the `hole', for some $i \neq h \in \{ 1, 2, \ldots, n \}$.  There are $n - 1$ options for this move, and we select the best one with respect to $\xi(\pi, h)$.  In addition to calculating $\xi(\pi, h)$, we also calculate $\phi(\pi)$ on every iteration to keep track of the best `full' solution found. A version of this ejection chain algorithm is presented in a preliminary report on this problem~\cite{our-chapter}.  For a formal description of our ejection chain algorithm see Algorithms~\ref{alg:ejection_chain} and \ref{alg:improvement}.

\begin{algorithm}[htb]
\caption{Ejection Chain Algorithm}
\label{alg:ejection_chain}

\begin{algorithmic}
\REQUIRE permutation $\pi^0$; maximum ejection chain length $\mathit{depth}$
\ENSURE Improved permutation with respect to $\phi()$

\STATE $\sigma \gets \pi^0$; $\pi \gets \pi^0$
\STATE $c \gets n$;  \COMMENT{\emph{$c$ counts non-improving iterations}}
\STATE $h \gets 1$
\WHILE {$c > 0$}
	\IF {$\mathit{improvement}(\pi, h, \sigma, \mathit{depth}) = 1$}
		\STATE $c \gets n$
		\STATE $\sigma \gets \pi$ as changed by $\mathit{improvement}(\pi, h, \sigma, \mathit{depth})$
	\ELSE
		\STATE $c \gets c - 1$
	\ENDIF
	\STATE \textbf{if} $h = n$ \textbf{then} $h \gets 1$ \textbf{else} $h \gets h + 1$; \COMMENT {\emph{Next hole position}}
\ENDWHILE
\RETURN $\sigma$
\end{algorithmic}
\end{algorithm}

\begin{algorithm}[htb]
\caption{Recursive Build of the Ejection Chain: $\mathit{improvement}(\pi, h, \sigma, d)$}
\label{alg:improvement}

\begin{algorithmic}
\REQUIRE permutation $\pi$; hole position $h$; the best solution found so far $\sigma$ (also an output parameter); the remaining search depth $d$
\ENSURE 1 --- improving ejection chain exists; 0 --- otherwise

\IF {$d = 0$}
	\RETURN 0
\ENDIF
\STATE $\xi_\text{p} \gets \xi(\pi, h)$
\STATE $j \gets 0$
\FOR {$i \gets 1, 2, \ldots, h - 1, h + 1, h + 2, \ldots, n$}
	\STATE $\pi' \gets \mathit{Swap}(\pi, h, i)$
	\IF {$\phi(\pi') > \phi(\sigma)$}
		\STATE $\sigma \gets \pi'$
		\RETURN 1
	\ENDIF
	\IF {$\xi(\pi', i) > \xi_\text{p}$}
		\STATE $j \gets i$
		\STATE $\xi_\text{p} \gets \xi(\pi', i)$
	\ENDIF
\ENDFOR

\IF {$j > 0$}
	\STATE $\pi' \gets \mathit{Swap}(\pi, h, j)$
	\IF {$\mathit{improvement}(\pi', j, \sigma, d - 1) = 1$}
		\RETURN 1
	\ENDIF
\ENDIF
\RETURN 0
\end{algorithmic}
\end{algorithm}

\subsection{Simulated Annealing}

Simulated Annealing (SA) is a stochastic optimisation technique widely used in the literature for solving various optimisation problems.  SA algorithm is similar to the randomised local search with the exception that the worsening moves can also be accepted.  

We implemented the standard SA based on the $N_{swap}(\pi)$ neighbourhood (see Algorithm~\ref{alg:sa}).  In each iteration, we swap two randomly selected elements in $\pi$.  If the obtained solution is better than $\pi$, we replace $\pi$ with that solution.  Otherwise, the probability of accepting the solution is $e^\frac{\phi(\pi') - \phi(\pi)}{T}$, where $T$ is the current temperature.  In each iteration, the temperature decreases linearly from a given initial value $T_0$ to 0.

\begin{algorithm}[htb]
\caption{Simulated Annealing improvement heuristic}
\label{alg:sa}

\begin{algorithmic}
\REQUIRE initial permutation $\pi^0$; initial temperature $T_0$; given time $\mathcal{T}$
\ENSURE Improved permutation with respect to $\phi()$

\STATE $\pi \gets \pi^0$
\WHILE {$\text{``elapsed time''} < \mathcal{T}$}
	\STATE Select $i \neq j \in \pi$ randomly with uniform distribution
	\STATE $\pi' \gets \mathit{Swap}(\pi, i, j)$
	\IF {$\phi(\pi') > \phi(\pi)$}
		\STATE $\pi \gets \pi'$
	\ELSE
		\STATE $T \gets T_0 \cdot \frac{\mathcal{T} - \mathit{``elapsed\;time''}}{\mathcal{T}}$
		\STATE $p \gets e^\frac{\phi(\pi') - \phi(\pi)}{T}$
		\STATE $r \gets$ random number uniformly distributed in $[0, 1]$
		\IF {$r < p$}
			\STATE $\pi \gets \pi'$
		\ENDIF
	\ENDIF
\ENDWHILE
\RETURN $\pi$
\end{algorithmic}
\end{algorithm}

\subsection{Tabu Search}

The Tabu Search (TS) meta-heuristic is a neighbourhood-based search methodology with a tabu list mechanism for escaping local maxima.  By storing certain features of the recent solutions in a \emph{tabu list}, TS avoids re-exploring previously visited areas of the search space, which, in turn, allows the algorithm to accept worsening solutions when it is at a local maximum.  A high-level description of the TS procedure is given in Algorithm~\ref{alg:tabu}.


\begin{algorithm}[htb]
\caption{Tabu Search improvement heuristic}
\label{alg:tabu}

\begin{algorithmic}
\REQUIRE initial permutation $\pi^0$; tabu list length $L$; given time $\mathcal{T}$
\ENSURE Improved permutation with respect to $\phi()$

\STATE Initialise an empty FIFO list $\mathcal{L}$
\STATE $\pi_\text{cur} \gets \pi^0$;  \COMMENT{$\pi_\text{cur}$ is the solution whose neighbourhood is explored}
\STATE $\pi^* \gets \pi$;  \COMMENT{$\pi^*$ is the best solution found in the neighbourhood of $\pi_\text{cur}$}
\WHILE {$\text{``elapsed time''} < \mathcal{T}$}
	\FORALL {$\{ i, j \} \subset \pi_\text{cur}$}
		\STATE $\pi' \gets \mathcal{Swap}(\pi_\text{cur}, i, j)$
		\IF {$\phi(\pi') > \phi(\pi)$; \COMMENT{The best known solution is improved}}
			\STATE $\pi \gets \pi'$; \COMMENT{Record as the best found solution}
			\STATE $\pi^* \gets \pi'$; \COMMENT{Update $\pi^*$ ignoring the possible tabu}
		\ELSIF {$\pi' \notin \mathcal{L}$ and $\phi(\pi') > \phi(\pi^*)$}
			\STATE $\pi^* \gets \pi'$
		\ENDIF
	\ENDFOR
	
	\STATE $\mathit{Insert}(\pi^*, \mathcal{L})$
	\IF {$|\mathcal{L}| > L$}
		\STATE $\mathit{Remove}(\mathcal{L})$
	\ENDIF
	\STATE $\pi_\text{cur} \gets \pi^*$ \COMMENT{Move to the best found solution}
\ENDWHILE
\RETURN $\pi$
\end{algorithmic}
\end{algorithm}

The efficiency of the TS algorithm significantly depends on the features to be kept in the tabu list.  For problems with permutation-based solution representation, it is a common practice to use pairs of recently modified elements and their positions as such features.  Any solution that has the saved elements at exactly the same positions is excluded from exploration in the next few iterations (tabu tenure) of the search.

However, our experiments have shown that such a TS implementation performs poorly on the DRPP\@.  For the explanation, observe that one SIDSP solution can be represented by many distinct DRPP solutions.  For instance, if a request $j$ is scheduled at its release time $r_j$, advancing $j$ in the permutation $\pi$ does not change the resulting schedule $S_\pi$.  Hence, simple constraints on the permutation $\pi$ do no guarantee that the SIDSP solution $S_\pi$ is excluded from the search, which might affect the ability of the TS to escape the local maximum.  In other words, tabu lists based on request positions work well in the space of permutations but fail in the space of schedules.

To make sure that the search does not return to the recently explored region of SIDSP solutions, we save features of schedules (rather than permutations) in the tabu list.  Each element of our tabu list includes the objective value and the average request tardiness $t(S) = \frac{1}{|S|} \cdot \sum_{j \in S} (S_j - r_j)$ of a recently explored solution $S = S_\pi$.  If both the objective value and the average tardiness of a new solution $S'$ are close to the ones in the list, such a solution is excluded from the search.  More formally, with respect to a tabu list element $(f, t)$, a solution $S'$ is tabu if (a) $\frac{|f(S') - f|}{f} \le \epsilon$ and (b) $\frac{|t(S') - t|}{t} \le \epsilon$, where $0 < \epsilon \ll 1$ is a tolerance parameter of the algorithm.

\section{Schedule Generator Algorithm}
\label{sec:ScheduleGenerator}

This section describes the schedule generation algorithm we use to produce a schedule from a given ordered subset of requests $V$.  Let us start by introducing some terminology and notations to simplify the discussion.

\subsection{Interval Sets}

Let $I$ be a finite set of non-intersecting intervals.  We call such a structure \emph{interval set}.  Since the elements of $I$ are non-intersecting intervals, $I$ can also be viewed as an ordered set with the natural order produced by the position of these intervals on the real line.  Thus, $I$ is represented as $I = \{[\ell_1,u_1], [\ell_2,u_2], \ldots, [\ell_v, u_v]\}$, where $\ell_1 < u_1 < \ell_2 < u_2 < \cdots < \ell_v < u_v$ and $v = |I|$.  We refer to the $k$-th interval in $I$ as $I_k$.

Consider two arbitrary intervals $[a,b]$ and $[\ell, u]$.  Let us introduce the \emph{subtraction} operation $[\ell, u] \circleddash [a, b]$ as follows:
$$
[\ell, u] \circleddash [a, b] =
\begin{cases}
\emptyset & \text{if } a \le \ell \text{ and } b \ge u, \\
\{ [\ell, a] \comma [b, u] \} & \text{if } a > \ell \text{ and } b < u, \\
\{ [b, u] \} & \text{if } a \le \ell \text{ and } b < u, \\
\{ [\ell, a] \} & \text{if } a > \ell \text{ and } b \ge u.
\end{cases}
$$
Now we can define subtraction $I \circleddash [a, b]$ for an interval $[a, b]$ and an interval set $I = \{[\ell_1,u_1] \comma [\ell_2,u_2] \comma \ldots \comma [\ell_v, u_v]\}$:
$$
I \circleddash [a, b] = \bigcup_{[\ell, u] \in I} [\ell, u] \circleddash [a, b] \,.
$$
Informally, one can think of the subtraction operation as a set subtraction where the interval set $I$ and the interval $[a, b]$ are represented as sets of points, and where the isolated points are excluded from the result.

An interval $[a,b]$ is said to be a \emph{subinterval} of the interval set $I$ if there exists $k$ such that $\ell_k \leq a < b \leq u_k$ and $[\ell_k,u_k] \in I$. This relationship is denoted by $[a,b] \Subset I$.

\subsection{Implementation of the Schedule Generator}

Let $V^*$ be an ordered subset of $V$.  Given $V^*$, we now present a \emph{schedule generation algorithm} to schedule requests in $V^*$ following the order prescribed in $V^*$.  The algorithm maintains several indicator interval sets representing channel availability at ground stations and antenna availability on the satellite.  After scheduling a request, the algorithm updates these indicator sets.  A high-level pseudo-code of the schedule generation algorithm is presented in Algorithm~\ref{alg:greedy}.  In what follows, we describe the details of each step of the algorithm.

\begin{algorithm}[htb]
\caption{The Schedule Generation Algorithm takes an ordered set of requests $V^* \subseteq V$ as an input and attempts to allocate each of the requests $j \in V^*$, following the prescribed order, to the earliest available time subject to all the operational constraints.  The output of the procedure is a feasible schedule.}
\label{alg:greedy}

\begin{algorithmic}
\REQUIRE An ordered set of requests $V^* \subseteq V$
\ENSURE A schedule $S$ as defined in Section~\ref{sec:problem}, i.e.\ a set of scheduled requests $S$ and associated downlink start times $S_j$ for each $j \in S$

\REPEAT
	\STATE Initialise the indicator sets, $S \gets \emptyset$, $V^{**} \gets V^*$ and $\mathit{restart} \gets 0$
	\WHILE {$V^{**} \neq \emptyset$ and $\mathit{restart} = 0$}
		\STATE Let $j$ be the first request in $V^{**}$
		\STATE Update $V^{**} \gets V^{**} \setminus \{ j \}$
		\STATE Find the earliest start time $x$ for request $j$
		\IF {no such $x$ exists}
			\IF {$\exists i \in V$ such that $(i, j) \in D$ or $(j, i) \in D$}
				\STATE $V^* \gets V^* \setminus \{ i, j \}$ and $V^{**} \gets V^{**} \setminus \{ i \}$
				\IF {$i \in S$}
					\STATE $\mathit{restart} \gets 1$
				\ENDIF
			\ENDIF
		\ELSE
			\STATE $S \gets S \cup \{ j \}$ and $S_j \gets x$
			\STATE Update the indicator sets
		\ENDIF
	\ENDWHILE
\UNTIL $\mathit{restart} = 0$

\STATE Assign the downlinks to the particular ground station channels and satellite antennas
\RETURN $S$
\end{algorithmic}
\end{algorithm}

As we mentioned above, the algorithm maintains several interval sets:
\begin{itemize}
	\item $A$ is an interval set representing the time intervals when both satellite antennas are available.
	\item $H$ and $F$ are interval sets indicating the time intervals when half-power and full-power downlinks can happen, respectively.
	\item $Q_g$ and $Q^1_g$ are interval sets indicating the availability of the ground station $g \in G$ in normal and high reliability visibility, respectively.  Note that we do not need the second pair of indicator sets for a two channel ground station.  Indeed, the number of channels in this case is not limiting since the number of simultaneous downlinks is constrained by the number of antennas.  Hence, resource availability of the two channel ground stations does not need to be tracked.
\end{itemize}

We initialise the indicator sets as follows: $A \gets H \gets F \gets \{[0, 24 \text{ hours}]\}$, and we set $Q_g$ to the normal visibility mask of $g$ and $Q^1_g$ to the high reliability visibility mask of $g$ for every ground station $g \in G$.  Note that for every $[\ell^1_k, u^1_k] \in Q^1_g$ there is an interval $[\ell_l, u_l] \in Q_g$ such that $\ell_l < \ell^1_k < u^1_k < u_l$.

If $j$ is to be downlinked to a half-power station $g \in G_1$, we use the following procedure to find the earliest time $x$ when it can be scheduled.  Let $Q \gets Q_g$ if $j$ requires normal reliability and $Q \gets Q^1_g$ otherwise.  Choose the smallest $k \in \{ 1, 2, \ldots, |H| \}$ and $l \in \{ 1, 2, \ldots, |Q| \}$ such that $|H_l \cap Q_k \cap [r_j, d_j]| \ge p_j$.  If no such $k$ and $l$ exist, request $j$ cannot be scheduled, and the algorithm proceeds to the next request.  Otherwise compute $[x, y] = H_k \cap Q_l \cap [r_j, d_j]$, schedule the request to tome $x$ and update the indicator sets as follows.  Let $X = [S_j, S_j + p_j]$.  Set $X \gets X \circleddash [\ell_k, u_k]$ for every $[\ell_k, u_k] \in A$.  Then $t \in X$ iff $t \in [S_j, S_j + p_j]$ and exactly one antenna was available at time $t$ before scheduling $j$.  Set $H \gets H \circleddash [\ell_k, u_k]$ for every $[\ell_k, u_k] \in X$ to reflect that the time intervals $X$ are no longer available to half-power downlinks.  Also set $A \gets A \circleddash [S_j, S_j + p_j]$ and $F \gets F \circleddash [S_j - \Delta, S_j + p_j + \Delta]$.  Finally, if $g$ is a one channel ground station, update $Q_g \gets Q_g \circleddash [S_j, S_j + p_j]$ and $Q^1_g \gets Q^1_g \circleddash [S_j, S_j + p_j]$ (recall that two channel ground stations are never limiting the number of simultaneous downlinks).

If $j$ is to be downlinked to a full-power station $g \in G_2$, we use another procedure to find the earliest time $x$ when it can be scheduled.  Let $Q \gets Q_g$ if $j$ requires normal reliability and $Q \gets Q^1_g$ otherwise.  Choose the smallest $k \in \{ 1, 2, \ldots, |F| \}$ and $l \in \{ 1, 2, \ldots, |Q| \}$ such that $|F_k \cap Q_l \cap [r_j, d_j]| \ge p_j$.  If no such $k$ and $l$ exist, request $j$ cannot scheduled, and the algorithm proceeds to the next request.  Otherwise compute $[x, y] = F_k \cap Q_l \cap [r_j, d_j]$, schedule the request to time $x$ and update the indicator sets as follows: $F \gets F \circleddash [S_j, S_j + p_j]$ and $H \gets H \circleddash [S_j - \Delta, S_j + p_j + \Delta]$.  Note that it is not necessary to update the indicator set $A$ since no downlink can happen if neither $H$ nor $F$ is available.

\bigskip

Our updating scheme of the interval sets $A$, $H$ and $F$ ensures that no antenna conflict arises and all the downlinks happen within the planning horizon.  By subtracting $[S_j - \Delta, S_j + p_j + \Delta]$ from $F$ for every half-power downlink request $j$, we guarantee that no full-power downlink can happen within $\Delta$ units of the downlink $j$.  Similarly, no half-power downlink can happen within $\Delta$ units of a full-power downlink $j$.  Also, the updating of the indicator sets $Q_g$ and $Q^1_g$ guarantees that no channel conflicts occurs and the downlinks obey visibility mask constraints.  The pre-processing of data assures that there is a gap of at least $\delta$ units between two consecutive downlinks.  Finally, dual requests constraint is satisfied as every time a conflict is detected (one of the requests is schedules while the other one cannot be scheduled), the procedure restarts.

Let us now analyse the complexity of the algorithm.  The primary operations in each iteration are: (1) to find the smallest $k$ and $l$ to satisfy certain condition, and (2) to update the indicator sets $A$, $H$, $F$, $Q_g$ and $Q^1_g$.  Note that the size of the indicator set $Q_g$ (and $Q^1_g$) for some $g \in G$ may increase but it is bounded by $O(n_g)$, where $n_g$ is the number of requests to be scheduled to the station $g$.  Similarly, the sizes of the indicator sets $A$, $H$ and $F$ are limited by $O(n)$.  Thus, operation (1) for these sets can be performed in $O(n)$ time by simultaneous scanning of the sets.  Operation (2) can also be performed in $O(n)$ time.  Indeed, we only need to update a fixed number of indicator sets, and updating each of them takes $O(n)$ time (for a downlink to a half-power setting ground station, manipulation with $X$ also requires only $O(n)$ time\footnote{Note that each of the $H \gets H \circleddash [\ell_k, u_k]$ for $[\ell_k, u_k] \in X$ operations needs only $O(1)$ time since all of these operations affect only one interval in $H$.}).  The number of iterations is at most $O(n^2)$ and, thus, the complexity of the algorithm  is $O(n^3)$.

In fact, the running time of this algorithm can potentially be reduced by smarter processing of dual requests.  In particular, one does not need to restart the algorithm every time a dual request constraint violation is observed.  It is enough to roll back the state of the algorithm to the point when one of the dual requests was scheduled.  However, that would not reduce the worst time complexity of the algorithm as the number of iterations would still be $O(n^2)$.  Indeed, each roll back requires $O(n)$ iterations as it may result in re-scheduling $O(n)$ requests, and there are $O(n)$ roll backs required in the general case.  Moreover, the roll back procedure would complicate the implementation (note that rolling back would require some form of restoration of the indicator sets).  Finally, we observed that the number of dual requests was low in our test instances, so we decided to restart the generation procedure for each detected dual request violation, as shown in Algorithm~\ref{alg:greedy}.

\section{Real-World Problem Instances}
\label{sec:testbed}

The algorithms presented in this paper have been tested on real data --- the RADARSAT-2 problem instances: (1) 10 \emph{low-density instances} (LD1---LD10) collected in Autumn 2011, each containing approximately 100 requests per planning horizon (24 hours); and (2) 10 \emph{high-density instances} (HD1---HD10) collected in August 2011, each containing approximately 300 requests.  There are at most ten ground stations involved in each instance.  Each ground station is visible between 4 to 10 times from the satellite during the planning horizon, depending on the ground station location.  Around 70\% of the ground stations have one channel and 30\% have two channels.

For the purpose of our experimental study, we were provided with the real downlink schedules implemented  for each of the low- and high-density instances.  The process currently in use for satellite mission planning includes two phases: (1) construction of the schedules with a priority rule-based algorithm\footnote{The details of that algorithm are unavailable to us.} and (2) human intervention.  The system operator modifies the machine-generated solutions with the aim of scheduling some additional downlink requests and satisfying additional considerations known to the operator at that time.  Such a solution may be imprecise.  For example, a human operator may sometimes use his/her judgement to schedule a downlink request even if a downlink goes beyond the prescribed visibility mask by a very small amount of time.  Such a solution would be infeasible as per our model as we use crisp visibility mask boundaries, as per satellite mission planning requirements.

We call such real schedules \emph{human-rescheduled} (H-R) and in the following section we compare them against our algorithms.

\section{Experimental Study}
\label{sec:results}

We implemented our algorithms in \verb!C++! and tested them on a PC with Intel Core i7-3820 CPU (3.60 GHz)\@.  The low-density and high-density instances discussed in Section~\ref{sec:testbed} were used in this experimental study.

Following an empirical parameter tuning procedure, we set the time given to each local search run within GRASP to 1 second, the EC search depth to 10, the simulated annealing initial temperature $T_0$ to $0.001$, the tabu list length $L$ to $4$, and the tabu tolerance parameter $\epsilon$ to $0.01$.

\subsection{Computational Results}
\label{sec:AlgComparison}

In this section we compare the performance of GRASP (named GR. in the tables below), Ejection Chain (EC), Simulated Annealing (SA) and Tabu Search (TS) algorithms (see Section~\ref{sec:dev_algorithm}) to the H-R schedules (see Section~\ref{sec:testbed}).  The computational results for low-density and high-density instances are reported in Tables~\ref{tab:comparison_low} and~\ref{tab:comparison_high}, respectively.  In our experiments, all our algorithms were given equal time for fair comparison.  Since EC is our only algorithm that does not have an explicit setting for the running time, we gave each of GRASP, SA and TS algorithms as much time as EC needed to terminate for each particular instance.

\begin{sidewaystable}
\centering
\footnotesize
\setlength{\tabcolsep}{0.5em}
\begin{tabular}{@{}l@{~~}r@{~~}r@{~~}r@{~~}c@{~~}r@{~~}r@{~~}r@{~~}r@{~~}r@{~~}c@{~~}r@{~~}r@{~~}r@{~~}r@{~~}r@{~~}c@{~~}r@{~~}r@{~~}r@{~~}r@{~~}r@{}}
\toprule

&&&&&\multicolumn{5}{@{}c@{}}{Unscheduled (urgent/total)}&&\multicolumn{5}{@{}c@{}}{Avg.\ tard.\ (urg.), sec}&&\multicolumn{5}{@{}c@{}}{Avg.\ tard., sec}\\
\cmidrule(r){6-10}
\cmidrule(r){12-16}
\cmidrule(){18-22}
Inst.&$|V|$&Urg.&Time&&H-R&GR.&EC&SA&TS&&H-R&GR.&EC&SA&TS&&H-R&GR.&EC&SA&TS\\
\cmidrule(){1-22}

LD1&110&31&1.7&&0 / 1&0 / 0&0 / 0&0 / 0&0 / 0&&0.3&0.0&0.0&0.0&0.0&&762&135&135&144&271\\
LD2&108&34&2.0&&0 / 1&0 / 0&0 / 0&0 / 0&0 / 0&&0.5&0.0&0.0&0.0&0.0&&1001&250&271&271&393\\
LD3&115&29&2.5&&0 / 1&0 / 0&0 / 0&0 / 0&0 / 0&&0.7&0.0&0.0&0.0&0.0&&881&95&94&98&219\\
LD4&105&25&2.2&&0 / 1&0 / 0&0 / 0&0 / 0&0 / 0&&0.7&0.0&0.0&0.0&0.0&&1613&296&457&296&680\\
LD5&104&32&0.7&&0 / 2&0 / 0&0 / 0&0 / 0&0 / 0&&0.4&0.0&0.0&0.0&0.0&&701&184&184&184&190\\
LD6&105&33&1.0&&0 / 1&0 / 0&0 / 0&0 / 0&0 / 0&&1.3&0.0&0.0&0.0&0.0&&202&119&119&119&135\\
LD7&90&30&0.7&&0 / 1&0 / 0&0 / 0&0 / 0&0 / 0&&3.1&0.0&0.0&0.0&0.0&&379&79&80&80&97\\
LD8&108&31&1.5&&0 / 1&0 / 0&0 / 0&0 / 0&0 / 0&&4.5&0.0&0.0&0.0&0.0&&234&103&105&102&118\\
LD9&110&28&1.8&&0 / 1&0 / 0&0 / 0&0 / 0&0 / 0&&11.3&0.0&0.0&0.0&0.0&&975&244&250&237&278\\
LD10&134&27&4.5&&0 / 0&0 / 0&0 / 0&0 / 0&0 / 0&&0.4&0.0&0.0&0.0&0.0&&564&65&65&66&199\\
\cmidrule(){1-22}

Average&108.9&30.0&1.8&&0.0 / 1.0&0.0 / 0.0&0.0 / 0.0&0.0 / 0.0&0.0 / 0.0&&2.3&0.0&0.0&0.0&0.0&&731.1&156.9&176.0&159.7&258.0\\

\bottomrule

\end{tabular}
\caption{Comparison of the heuristic algorithms for the low-density instances.}
\label{tab:comparison_low}

\bigskip
\bigskip

\begin{tabular}{@{}l@{~~}r@{~~}r@{~~}r@{~}c@{~}r@{~~}r@{~~}r@{~~}r@{~~}r@{~}c@{~}r@{~~}r@{~~}r@{~~}r@{~~}r@{~}c@{~}r@{~~}r@{~~}r@{~~}r@{~~}r@{}}
\toprule

&&&&&\multicolumn{5}{@{}c@{}}{Unscheduled (urgent/total)}&&\multicolumn{5}{@{}c@{}}{Avg.\ tard.\ (urg.), sec}&&\multicolumn{5}{@{}c@{}}{Avg.\ tard., sec}\\
\cmidrule(r){6-10}
\cmidrule(r){12-16}
\cmidrule(){18-22}
Inst.&$|V|$&Urg.&Time&&H-R&GR.&EC&SA&TS&&H-R&GR.&EC&SA&TS&&H-R&GR.&EC&SA&TS\\
\cmidrule(){1-22}

HD1&211&34&29.4&&0 / 53&0 / 19&0 / 19&0 / 18&0 / 24&&0.3&0.0&0.0&0.0&0.0&&2013&581&450&462&630\\
HD2&301&35&135.8&&0 / 123&0 / 67&0 / 65&0 / 64&0 / 65&&0.6&0.0&0.0&0.0&0.0&&3788&3087&2943&3074&5446\\
HD3&324&42&432.9&&0 / 120&0 / 51&0 / 46&0 / 47&0 / 52&&0.5&0.0&0.0&0.0&0.0&&2329&2693&2397&2183&2745\\
HD4&294&41&149.1&&2 / 116&0 / 57&0 / 54&0 / 51&0 / 61&&0.5&0.0&0.0&0.0&0.0&&2025&1529&1201&1326&1864\\
HD5&260&36&284.3&&0 / 85&0 / 50&0 / 49&0 / 48&0 / 49&&0.5&0.0&0.0&0.0&0.0&&2621&899&1013&1294&1390\\
HD6&356&46&468.8&&2 / 130&1 / 63&1 / 58&1 / 58&1 / 68&&0.4&0.0&0.0&0.0&0.0&&2555&1878&1702&1885&3921\\
HD7&259&38&124.6&&2 / 87&0 / 31&0 / 32&0 / 30&0 / 30&&0.4&0.0&0.0&0.0&0.0&&2407&1994&1658&1941&2557\\
HD8&304&54&229.4&&1 / 98&0 / 35&0 / 35&0 / 34&0 / 36&&0.4&0.0&0.0&0.0&0.0&&2230&1991&1396&1250&2345\\
HD9&278&39&134.9&&3 / 90&0 / 37&0 / 37&0 / 35&0 / 39&&0.5&0.5&0.5&0.5&0.5&&2463&1085&908&761&2304\\
HD10&230&25&83.4&&0 / 80&0 / 28&0 / 28&0 / 28&0 / 30&&0.5&0.0&0.0&0.0&0.0&&2880&2092&2081&2077&2342\\
\cmidrule(){1-22}

Avg.&281.7&39.0&207.3&&1.0 / 98.2&0.1 / 43.8&0.1 / 42.3&0.1 / 41.3&0.1 / 45.4&&0.5&0.1&0.1&0.1&0.1&&2531.0&1782.9&1574.8&1625.5&2554.4\\

\bottomrule

\end{tabular}

\caption{Comparison of the heuristic algorithms for the high-density instances.}
\label{tab:comparison_high}
\end{sidewaystable}

The columns of Tables~\ref{tab:comparison_low} and~\ref{tab:comparison_high} are as follows (from left to right): the instance name, the number of downlink requests $|V|$, the number of urgent requests, the time given to each of our algorithms, the number of unscheduled (urgent/total) requests for each of the algorithms, the average tardiness of the urgent requests for each of the algorithms, and the overall average tardiness for each of the algorithms.  The tardiness of a request $j$ is measured as $S_j - \tau_j$, where $\tau_j$ is the earliest time $j$ can be downlinked subject to no other downlinks are scheduled.

It follows from the results of our computational experiments that the low-density instances are relatively easy to solve.  Observe that each of our algorithms (GRASP, EC, SA and TS) scheduled all the requests and achieved 0.0~seconds urgent request tardiness for each instance.  In terms of the overall average tardiness, the GRASP and the SA algorithms are the winners.

The high-density instances are much harder to solve.  Each of the algorithms left several requests unscheduled for each of the instances.  In terms of urgent requests, all our algorithms performed very similarly.  In terms of overall performance, SA and EC are the leaders, followed by GRASP\@.

Our algorithms clearly outperformed the H-R solutions.  For the low-density instances, the H-R solutions left, on average, one unscheduled request while all of our algorithms managed to schedule all the requests.  For the high-density instances, the H-R solutions left, on average, 98.2 requests unscheduled, which is more than twice compared to any of our algorithms.  Our best algorithms also significantly decreased the average tardiness (both overall and for urgent requests) compared to the H-R solutions.

It is worth noting that our construction heuristic described in Section~\ref{sec:InitialSoln} also outperformed the H-R solutions.  For example, for the high-density instances, it left only 56.9 requests unscheduled, on average.  Compare it to the 98.2 unscheduled requests in the H-R solutions.  However, each of our meta-heuristics significantly improved the results of the construction heuristic, leaving only 41--45 requests unscheduled.

\subsection{Effect of Given Time Parameter}

In this section we test the effect of varying the running time of our algorithms.  The GRASP, SA and TS algorithms have explicit parameters to adjust their running time.  The EC algorithm has only one parameter, the search depth, that might affect the running time of EC\@.  In this experiment, we ran GRASP, SA and TS given $4, 8, 16, \ldots, 8192$~seconds and EC with $\mathit{depth} = 2, 5, 10, 20, 50$.  The results for each algorithm and each setting were averaged over all the high-density instances.

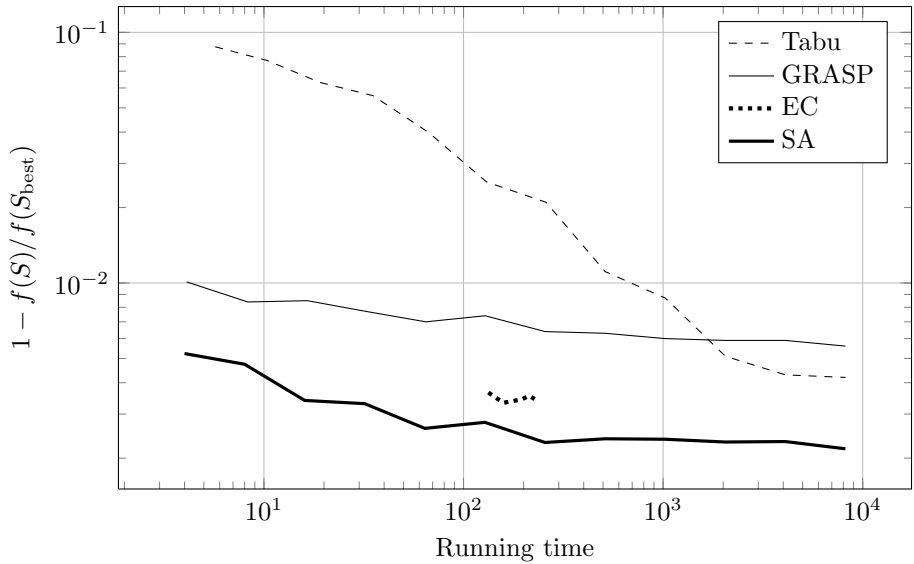
\begin{figure}
\begin{tikzpicture}
	\begin{loglogaxis}[
		width=\textwidth,
		height=8cm,
		legend pos=north east,
		legend cell align=left,
		xlabel={Running time},
		ylabel={$1 - f(S) / f(S_\text{best})$},
		grid=major
	]
	\addplot[dashed, black] coordinates {
		(5.7,   0.0876)
		(10.4,  0.0772)
		(19.0,  0.0632)
		(35.2,  0.0559)
		(66.8,  0.0399)
		(130.5, 0.0253)
		(258.7, 0.0210)
		(512.9, 0.0111)
		(1026.9,0.0087)
		(2052.0,0.0051)
		(4098.9,0.0043)
		(8194.7,0.0042)
	};
	\addlegendentry{Tabu}
	\addplot[black] coordinates {
		(4.1,   0.0101)
		(8.3,   0.0084)
		(16.5,  0.0085)
		(32.6,  0.0077)
		(64.7,  0.0070)
		(128.4, 0.0074)
		(256.3, 0.0064)
		(512.4, 0.0063)
		(1024.3,0.0060)
		(2048.6,0.0059)
		(4096.6,0.0059)
		(8192.5,0.0056)
	};
	\addlegendentry{GRASP}
	\addplot[ultra thick, dotted, black] coordinates {
		(133.1265734,0.00365620269203123)
		(158.3160549,0.0033265901637823)
		(195.5254133,0.00343071664208311)
		(213.3232317,0.00354767784769461)
		(227.1826964,0.00344127847787469)
	};
	\addlegendentry{EC}
	\addplot[very thick, black] coordinates {
		(4.0,   0.00523)
		(8.0,   0.00474)
		(16.0,  0.00340)
		(32.0,  0.00330)
		(64.0,  0.00263)
		(128.0, 0.00278)
		(256.0, 0.00231)
		(512.0, 0.00239)
		(1024.0,0.00238)
		(2048.0,0.00232)
		(4096.0,0.00233)
		(8192.0,0.00218)
	};
	\addlegendentry{SA}	
	\end{loglogaxis}
\end{tikzpicture}

\caption{Performance of the GRASP, EC, SA and TS algorithms on a range of settings (given time for GRASP, SA and TS, and search depth $\mathit{depth}$ for EC).  The results are averaged over 10 runs (one run for each of the high-density instances).  The vertical axis shows how far a solution is from the best known solution.}
\label{fig:comparison}
\end{figure}

The results of the experiment are reported in Figure~\ref{fig:comparison}.  An important observation is that neither the running time nor the solution quality of the EC algorithm notably depend on the value of parameter $\mathit{depth}$.  The other algorithms' quality significantly improves when given more time.  TS shows relatively poor performance when given little time which, however, rapidly improves with the increase of the running time.  GRASP and SA are less sensitive to the given time.  The winning algorithm is clearly SA: it outperforms other algorithms on the whole range of given times, producing reasonable solutions in just a few seconds and being able to gradually improve the solution quality.

\subsection{Effect of Parameter $\alpha$}

\begin{figure}[htb]
\centering
\begin{tikzpicture}
\begin{axis}[
width=1\textwidth,
height=7cm,
xlabel={Average number of unscheduled requests per day},
ylabel={Average delay per request, hours},
nodes near coords,
every node near coord/.style={anchor=west,font=\footnotesize},
enlargelimits=0.1,
grid=major,
]
\addplot[scatter,only marks,point meta=explicit symbolic]
coordinates {
(40.73,4.70538840837103) [0]
(40.65,4.43790586918714) [0.01]
(41.17,4.22346024804396) [0.2]
(41.64,4.05778774864539) [0.5]
(42.25,3.91919673075254) [0.7]
(43.14,3.81199913774697) [0.8]
(44.17,3.70787706076115) [0.9]
(44.77,3.64274895065522) [0.95]
(45.65,3.5593544875655) [1]
};
\end{axis}
\end{tikzpicture}
\caption{Demonstration of how the value of $\alpha$ influences the trade-off between average tardiness and the number of unscheduled requests.  The number on the right of each point is the value of $\alpha$.  Each point was obtained by averaging the results of 30-second runs of the SA algorithm for 10 different seed values and all the high-density instances.}
\label{fig:alpha}
\end{figure}
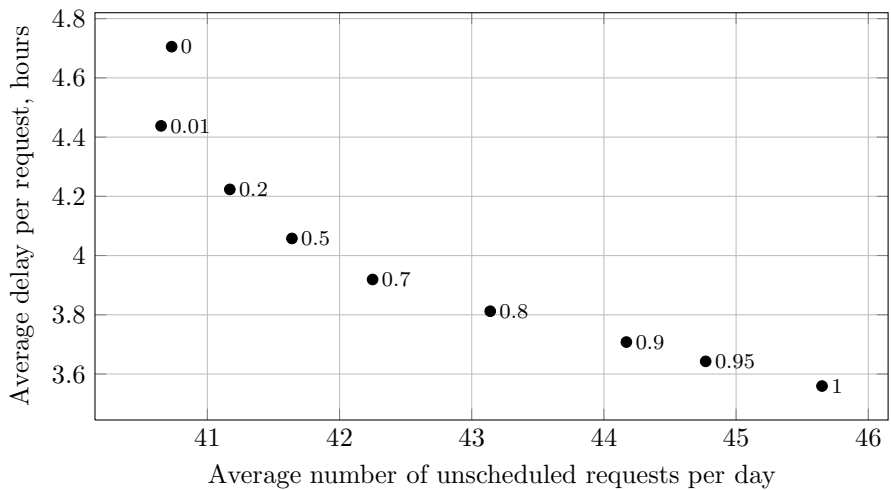

Recall that the objective function for the SIDSP uses a parameter $\alpha$ that represents the importance of tardiness minimisation compared to the importance of scheduling as many requests as possible.  The value of $\alpha$ is irrelevant as long as all the requests are scheduled (like in our solutions of the low-density instances).  However, if the problem is over-subscribed, the parameter $\alpha$ controls the trade-off between the number of scheduled requests and their tardiness.  Figure~\ref{fig:alpha} shows how the performance of the SA algorithm depends on $\alpha$.

As expected, increasing the value of $\alpha$ puts emphasis on tardiness minimisation at the cost of scheduling fewer requests.  Nevertheless, the range of solutions produced for different values of $\alpha$ is relatively small.  Hence, the results obtained in this paper for $\alpha = 0.5$ are expected to hold for different values of $\alpha$.

It is also interesting to note that the solutions obtained for $\alpha = 0$ are, on average, dominated by the solutions obtained for $\alpha = 0.01$.  We link it to the changing fitness landscape.  Observe that, for $\alpha = 0$, the objective value of a solution depends only on the set $S$ but not the times $S_j$.  Hence, many distinct solutions sharing the same set $S$ cannot be distinguished by the solver.  That creates the so-called high neutrality of the fitness landscape that is known to reduce the performance of optimisation heuristics~\cite{marmion}.

%

\section{Conclusions}
\label{sec:conclusion}

In this paper, we formalised the satellite downlink scheduling problem and proposed a flexible solution approach separating the details of the problem-specific constraints from the optimisation mechanism.  That significantly simplified the design and implementation of optimisation meta-heuristics.  We tested several standard search techniques, including GRASP, Ejection Chain, Simulated Annealing and Tabu Search, and chose Simulated Annealing as the most efficient algorithm.


Our heuristic achieved very promising results comparing to the solutions obtained by the currently implemented method.  For situations and problem instances where downlink scheduling was oversubscribed, the number of unscheduled downlink requests was halved.  Although many of unsuccessful downlinks could have been `background' acquisitions of little importance, the increase in the downlink throughput was significant and could be used for generating additional imagery products, increasing the number of customers, and thus increasing the company's profit.  Moreover, our algorithms produced solutions with lower total tardiness, on average.



This research shows how a complicated model can be combined with modern optimisation heuristic methods in such a way that the resulting system is easy to use and maintain.  The proposed approach gives an opportunity to handle other variations of the problem as well as potential additional constraints without re-implementing the optimisation part, which is a vital requirement of industrial systems.

\section*{Acknowledgement} We are thankful to the referees and the associate editor for helpful comments which improved the presentation of the paper. This work was supported by the NSERC CRD grant 411294-2010 awarded to Abraham P.~Punnen with the MDA Systems Ltd.\ as the collaborating industrial partner.


\end{document}